\documentclass[11pt]{article}

\usepackage{graphicx}
\usepackage{epsfig}
\usepackage{epstopdf}

\usepackage{amssymb, amsmath, amsthm}
\usepackage{float}

\newtheorem{dfn}{Definition}
\newtheorem{defn}[dfn]{Definition}

\newtheorem{rem}[dfn]{Remark}

\newtheorem{thm}{Theorem}
\newtheorem{lem}{Lemma}
\newtheorem{cor}{Corollary}

\begin{document}

\title{Flipping Heegaard splittings of Seifert fibered spaces}

\author{Jennifer Schultens}

\maketitle

\begin{abstract}
A Heegaard splitting of a $3$-manifold is flippable if there is an isotopy that interchanges the two sides of the Heegaard splitting.  We explore which Heegaard splittings of Seifert fibered spaces are flippable. 
\end{abstract}

\section{Overview}

Heegaard splittings of $3$-manifolds are often studied up to isotopy of their splittings surfaces, a point of view that serves to highlight the positioning of the Heegaard surface in a $3$-manifold.  Of further interest, however, is to understand whether or not the two sides of a Heegaard surface, the two handlebodies into which the Heegaard surface partitions the $3$-manifold, can be interchanged via isotopy.  We consider the case of Heegaard splittings of Seifert fibered spaces, some of which, it turns out, can be flipped while others cannot.  Whether or not a Heegaard splitting is flippable depends on the existence of a certain ${\mathbb Z}_2$ subgroup of the Goeritz group.  We here exploit our understanding of Heegaard splittings of Seifert fibered spaces along with the obstruction described in \cite{johnson2010sg} to give a complete classification of flippability of Heegaard splittings of Seifert fibered spaces.  

\section{Preliminaries}

In addition to providing the key definitions for considering flippability of Heegaard splittings,  we will also provide motivating examples.  Henceforward, we only consider orientable $3$-manifolds.  

\begin{defn}
A {\em handlebody} is a $3$-manifold homeomorphic to a $3$-dimensional regular neighborhood of a graph.
\end{defn}

\begin{defn}
A {\em Heegaard splitting} of a $3$-manifold $M$ is a decomposition, $M = H_1 \cup_S H_2$ into two handlebodies such
that $H_1 \cap H_2 = \partial H_1 = \partial H_2$ is a closed surface $S$ called the {\em splitting surface}.  Two Heegaard splittings are considered {\em equivalent}, if their splittings surfaces are isotopic.  
\end{defn}

\noindent
{\bf Example 1:}  The $3$-sphere (thought of as the unit sphere in ${\mathbb R}^4$) has an upper hemisphere, ${\mathbb B}_n$, and a lower hemisphere, ${\mathbb B}_s$.  Each hemisphere is a $3$-ball, {\it i.e.,} a regular neighborhood of a point, hence a handlebody.  Moreover, ${\mathbb B}_n \cap {\mathbb B}_s = {\mathbb S}^2$.  Thus ${\mathbb S}^3 = {\mathbb B}_n \cup_{{\mathbb S}^2} {\mathbb B}_s$ is a Heegaard splitting.  

\begin{lem} \label{genus0}
The genus $0$ Heegaard splitting of ${\mathbb S}^3$ is unique.  
\end{lem}

\proof
Every $3$-ball is a regular neighborhood of a point.  Any point can be isotoped into any other point.  By the isotopy uniqueness of regular neighborhoods, all genus $0$ Heegaard splittings of ${\mathbb S}^3$ are therefore isotopic.  
\qed

\noindent
{\bf Example 2:}  The torus ${\mathbb T} = {\mathbb S}^1 \times {\mathbb S}^1$ divides the $3$-sphere (thought of as the unit sphere in ${\mathbb R}^4$) into two solid tori $V_{io}$ and  $V_{oi}$.  Moreover, $V_{io} \cap V_{oi} = {\mathbb T}^2$.  Thus ${\mathbb S}^3 = V_{io} \cup_{{\mathbb T}^2} V_{oi}$ is a Heegaard splitting.  

\begin{lem} \label{genus1}
The genus $1$ Heegaard splitting of ${\mathbb S}^3$ is unique.  
\end{lem}

\proof
Every solid torus is a regular neighborhood of its core, a simple closed curve.  It is not hard to show that a simple closed curve in ${\mathbb S}^3$ is the core of a solid torus in a genus $1$ Heegaard splitting of ${\mathbb S}^3$ if and only if it is unknotted.  Since all unknots are isotopic, the isotopy uniqueness of regular neighborhoods also establishes that all genus $1$ Heegaard splittings of ${\mathbb S}^3$ are isotopic.  
\qed

\begin{defn}
A {\em stabilization} of a Heegaard splitting $M = H_1 \cup_S H_2$ is the Heegaard splitting whose splitting surface is defined by the pairwise connected sum $(M, S) \# ({\mathbb S}^3, T) \# \dots \# ({\mathbb S}^3, T)$ with a nonzero number of $({\mathbb S}^3, T)$ summands.  
\end{defn}

\begin{defn}
A Heegaard splitting $M = H_1 \cup_S H_2$ is said to be {\em flippable} if there is an isotopy of M that takes $S$ to itself but interchanges $H_1$ and $H_2$.  (This is equivalent to requiring that there is an isotopy that takes the oriented surface $S$ to itself but with the opposite orientation.)
\end{defn}

The proof of Lemma \ref{genus0} also establishes the following corollary:

\begin{cor} 
The genus $0$ Heegaard splitting of ${\mathbb S}^3$ is flippable.  
\end{cor}

The proof of Lemma \ref{genus1} also establishes the following corollary:

\begin{cor} 
The genus $1$ Heegaard splitting of ${\mathbb S}^3$ is flippable.  
\end{cor}

The following theorem is easy to deduce:

\begin{thm} \label{flipstable}
Stabilizations of flippable Heegaard splittings are flippable.  
\end{thm}

\section{Motivating example}

\begin{figure}[h]
\includegraphics[scale=.4]{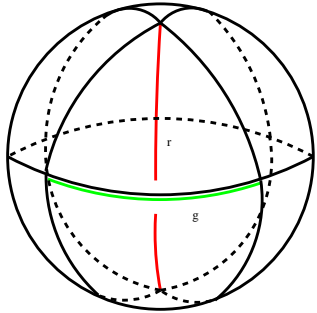}
\centering
\caption{\sl The lens space $L(4,1)$}
\label{lens}
\end{figure}

Lens spaces provide an interesting playground for flipping Heegaard splittings.  We obtain $L(p,q)$ by identifying appropriate sectors of the upper and lower hemisphere of ${\mathbb S}^2$ viewed as the boundary of the $3$-ball.  See Figure \ref{lens}.  Bonahon and Otal proved that Heegaard splittings of lens spaces are standard. See \cite{BO1982}.  In Figure \ref{lens}, the cores of the solid tori in the genus $1$ Heegaard splitting are labeled $r$ (red) and $g$ (green).  

\begin{figure}[h]
\includegraphics[scale=.4]{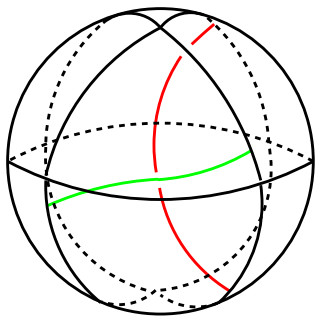}
\centering
\caption{\sl Isotoping cores of handlebodies in $L(4,1)$}
\label{lens1}
\end{figure}

\begin{thm} \label{lpq}
The genus $1$ Heegaard splitting of $L(p, q)$ is flippable if and only if $q = 1 \; (mod\;p)$.  
\end{thm}

\proof
Flipping the genus $1$ Heegaard splitting of $L(p, q)$ amounts to interchanging the two solid tori in the Heegaard splitting.  Since each solid torus is a regular neighborhood of its core, this is possible if and only if the two cores can be interchanged via an isotopy. 
Figure \ref{lens1} shows the beginning of an isotopy of $r$ and $g$.   Given the identifications of the sectors in $L(p,q)$, we see that $r$ is isotopic to $q$ times $g$.  When $q = 1\;(mod\;p)$, the isotopy transforms $r$ into $g$.  Likewise, $g$ can be transformed into $r$ via an isotopy and with a bit of care the two isotopies can be performed without obstructing each other.  

On the other hand, since $r$ is isotopic to a curve that wraps around $g$ exactly $q$ times, the two curves are not isotopic when $q \neq 1\;(mod\;p)$.  Thus for $q \neq 1\;(mod\;p)$, the genus $1$ Heegaard splitting is not flippable.  
\qed

\begin{rem}
The lens space $L(2, 1)$ is also known as ${\mathbb R}P^3$.  By Theorem \ref{lpq}, its genus $1$ Heegaard splitting is flippable.  
\end{rem}

\section{Flipping Heegaard splittings of products}

Heegaard splittings of product manifolds were studied in \cite{BO1982} and \cite{Schultens1993}.  

\begin{figure}[h]
\includegraphics[scale=.6]{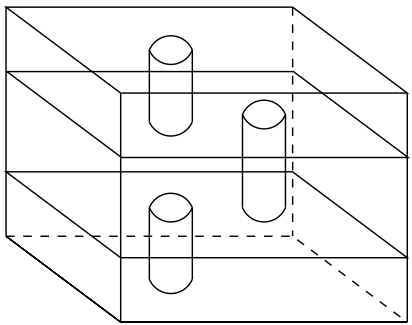}
\centering
\caption{\sl The standard Heegaard splitting of the $3$-torus}
\label{t3hs}
\end{figure}

\begin{defn} \label{prod}
Let $S$ be a closed connected orientable surface and set $M = S \times {\mathbb S}^1$.  The {\em standard} Heegaard splitting of $M$ is obtained by partitioning ${\mathbb S}^1$ into two subintervals $I_1 \cup I_2$ that intersect only in their endpoints; choosing two distinct points $p_1, p_2 \in S$ with disjoint open regular neighborhoods $\eta(p_1), \eta(p_2) \subset S$; and setting $$H_1 = (S - \eta(p_2)) \cup (N_1 \times {\mathbb S}^1)$$ $$H_2 = (S - \eta(p_1)) \cup (N_2 \times {\mathbb S}^1)$$
Where $N_1 = closure(\eta(p_1))$ and $N_2 = closure(\eta(p_2))$.  
\end{defn}

\begin{rem} \label{products1}
The Heegaard splitting described in Definition \ref{products} is isotopic to the Heegaard splitting obtained from the standard cell decomposition of $M$ by taking $H_1$ to be a regular neighborhood of the $1$-skeleton and $H_2$ to be a (the appropriate) regular neighborhood of the $1$-skeleton of the dual decomposition.  
\end{rem}

\begin{thm} (Boileau-Otal, Schultens) \label{products}
All Heegaard splittings of $S \times {\mathbb S}^1$ are standard.
\end{thm}

In \cite{BO1982}, Boileau and Otal showed that irreducible Heegaard splittings of the $3$-torus are isotopic to the standard Heegaard splitting and their theorem was extended to all $3$-manifolds of the form $(closed \; orientable \; surface) \times {\mathbb S}^1$ by the author of this paper.  See \cite{Schultens1993}.

\begin{thm} \label{prodflip}
Heegaard splittings of product manifolds are flippable.
\end{thm}

\proof
To interchange $H_1$ and $H_2,$ let $\alpha$ be a simple path in $S$ from $p_1$ to $p_2$ and let $\beta$ be a simple path in $S$ from $p_2$ to $p_1$ that is disjoint from $\alpha$ except at its endpoints.  An isotopy that rotates $S \times \partial I_1$ halfway around the ${\mathbb S}^1$ factor interchanges $I_1$ and $I_2$.  This isotopy can be accomplished while isotoping $p_1$ to $p_2$ along $\alpha$ and $p_2$ to $p_1$ along $\beta$.   The isotopy described then interchanges $H_1$ and $H_2$.

By Theorem \ref{products}, all Heegaard splittings of $(closed \; orientable \; surface) \times {\mathbb S}^1$ are stabilizations of the standard Heegaard splitting.  By Theorem \ref{flipstable}, stabilizations of flippable Heegaard splittings are flippable.  Therefore all Heegaard splittings of these manifolds are flippable.  
\qed

\section{Heegaard splittings of Seifert fibered spaces}

Seifert fibered spaces where originally identified and described by H. Seifert in \cite{SeifertBook}.  See also W. Heil's translation, \cite{HeilTrans}.  We describe their basic building blocks and state key results.

\begin{defn} \label{FST}
A {\em fibered solid torus} is a compact $3$-manifold obtained from a solid cylinder  ${\mathbb D}^2 \times [0,1]$ by glueing ${\mathbb D} \times \{0\}$ 
to ${\mathbb D} \times \{1\}$ after a rotation by a rational multiple of $2\pi.$  More specifically, the rotation is by $ \frac{2\pi \nu}{\mu}$, where $\mu, \nu \in \mathbb{Z}$ and $g.c.d.(\mu,\nu) = 1$.  We denote the resulting solid torus by $V(\nu,\mu).$  By convention, we require that $0 < \nu < \mu.$  
A simple closed curve resulting by identification of the endpoints of intervals of the form $\{y\} \times [0, 1], y \in {\mathbb D}^2$ is called a {\em fiber}.

An {\em exceptional} fiber is a fiber in $V(\nu,\mu)$, for $\mu > 1$, that results from identifying $\{(0,0)\} \in {\mathbb D}^2 \times [0,1]$ to $\{(0,1)\} \in {\mathbb D}^2 \times [0,1]$.  All other fibers are {\em regular} fibers. 
\end{defn}

\begin{figure}[h]
\includegraphics[scale=.4]{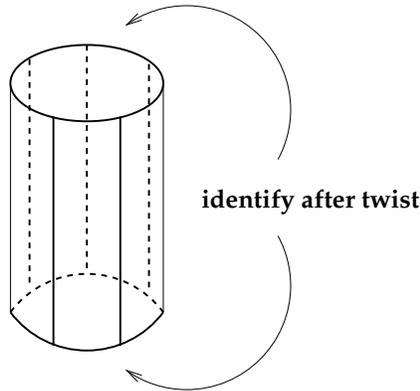}
\centering
\caption{\sl A fibered solid torus}
\label{fibtor}
\end{figure}

\begin{defn}\label{DSFS}
A {\em Seifert fibered space} $M$ is a compact connected $3$-manifold that admits a decomposition into
disjoint simple closed curves each of which has a neighborhood that is 
homeomorphic to a fibered solid torus via a homeomorphism that takes simple closed curves to fibers.  The simple closed curves into which $M$ decomposes are called {\em fibers} of $M.$ A particular decomposition of $M$ into fibers is called a {\em fiberation} of $M.$  A fiber of $M$ is {\em exceptional} if it is an exceptional fiber in a fibered solid torus neighborhood and {\em regular} otherwise.  See Figure \ref{fibtor}.  
\end{defn}

\begin{defn} \label{orbit}
Given a Seifert fibered space $M$, we form a quotient space, $Q$, by identifying each fiber to a point.   The quotient map is denoted by $p: M \rightarrow Q.$  Topologically, the quotient space is an orbifold called the {\em base orbifold}.  Topologically, $Q$ is a surface.  However, if nearby regular fibers wrap around the exceptional fiber $e$ exactly $\mu$ times, then we declare $p(e)$ to be a {\em cone point} of $Q$ of {\em multiplicity} $\mu.$

If the underlying surface of $Q$ is a sphere, then we say that $M$ has a {\em spherical} base orbifold.  
\end{defn}

\begin{figure}[h]
\includegraphics[scale=.5]{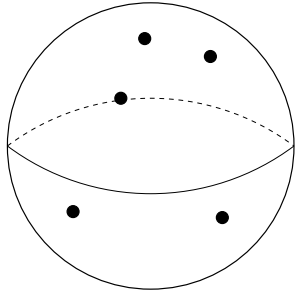}
\centering
\caption{\sl A spherical base orbifold with $4$ cone points}
\label{BaseOrbifold}
\end{figure}

\begin{defn}
A Seifert fibered space $p: M \rightarrow Q$ is {\em totally} orientable if $M$ and $Q$ are orientable.
\end{defn}

\begin{rem} \label{to}
Henceforward we will consider only totally orientable Seifert fibered spaces.  
\end{rem}

The following theorem describes how incompressible surfaces can be positioned with respect to the fibers of a Seifert fibered spaces.  

\begin{thm}\cite[VI.34]{JacoBook}   \label{Jaco}
(Jaco) Let $F$ be a connected, two-sided, essential surface in an orientable Seifert fibered space $M.$ Then one of the following holds: 

(i) $F$ is non separating in $M$ and is a fiber in a fiberation of $M$ as a surface bundle over the circle;

(ii) $F$ separates $M$ and $M = M_1 \cup M_2,$ where $\partial M_i = F$ and $M_i$ is a twisted $I$-bundle over a compact surface;

(iii) $F$ is an annulus or torus and $F$ is saturated, {\em i.e.,} consists of fibers, in some Seifert fiberation of $M.$
 
 \end{thm}
 
 \begin{defn}
 Let $M$ be a Seifert fibered space.  A surface $S$ in $M$ is {\em horizontal} if it is everywhere transverse to the fiberation.   It is {\em vertical} if it consists of fibers.  \end{defn}

The words ``horizontal" and ``vertical" are also used in the context of Heegaard splittings, but the descriptions are more complicated.  It should be noted that horizontal Heegaard splittings exist only in exceptional circumstances.  

\begin{defn} \label{hor}
Let $M$ be a Seifert fibered space.  A Heegaard splitting $M = H_1 \cup_S H_2$ is {\em horizontal} if, after isotopy, if necessary, there is a fiber $f$ of $M$ lying in $S$ such that $S \cap (M - \eta(f))$ is horizontal.  
\end{defn}

\begin{rem} \label{horiHS}
Horizontal Heegaard splittings can be defined more concretely:  Given that $S \cap (M - \eta(f))$ is horizontal and has two boundary components, our assumption that $M$ is totally orientable (see Remark \ref{to}) together with Theorem \ref{Jaco} implies that, after isotopy, if necessary, $S \cap (M - \eta(f))$ is a mapping torus over a surface $\hat S$, {\it i.e.,} $$M - \eta(f) =  (\hat{S} \times [-1,1])/\sim$$ and $S \cap (M - \eta(f))$ consists of two copies of $\hat S$.  Furthermore, $S$ intersects the regular neighborhood of $f$ in an annulus that is a bicollar of $f$.  
\end{rem}

\begin{figure}[ht]
\includegraphics[scale=.5]{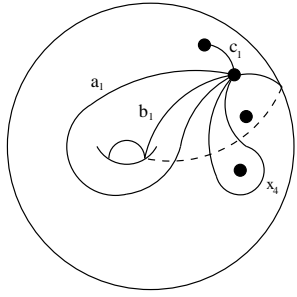}
\centering
\caption{\sl Arcs for a vertical Heegaard splitting}
\label{vertical}
\end{figure}

\begin{defn}
Let $M$ be a Seifert fibered space with $l$ exceptional fibers.  A {\em vertical} Heegaard splitting is a Heegaard splitting that, after isotopy, arises from the following construction:  
\begin{enumerate}
\item Case 1 ($M$ is a product) In this case the construction is as in Definition \ref{prod}, see also Remark \ref{products1}.  
\item Case 2 ($l \geq 2$) Denote the exceptional fibers of $M$ by $e_1, \dots, e_l$.  Choose a proper subset $e_1, \dots, e_n$.  In the base orbifold, choose a collection of arcs $ c_1, \dots, c_{n-1}, x_n, \dots, x_l, a_1, b_1, \dots, a_g, b_g$, see Figure \ref{vertical}.  Abusing notation, we continue to denote lifts of these arcs to $M$ by the same symbols.  Set $$\Sigma_1 = e_1 \cup \dots \cup e_n \cup c_1 \cup \dots \cup c_{n-1} \cup x_{n+2} \cup \dots \cup x_l \cup a_1 \cup b_1 \cup \dots \cup a_g \cup b_g$$ and take $H_1 = N(\Sigma_1)$.  Take $H_2$ to be the closure of the complement of $H_1$.  Alternatively, set $$\Sigma_2 =  e_l \cup \dots \cup e_{n+2} \cup c_{n+1}' \cup \dots \cup c_{l-1}' \cup x_1' \cup \dots \cup x_{n-1}' \cup a_1' \cup b_1' \cup \dots \cup a_g' \cup b_g'$$ and $H_2 = N(\Sigma_2)$ for an appropriately chosen alternate collection of arcs $ c_{n+1}', \dots, c_{l-1}', x_{1}', \dots, x_{n-1}', a_1', b_1', \dots, a_g', b_g'$ based at $p(e_n)$.    
\item Case 3 ($M$ is a nontrivial circle bundle or $l = 1$) Proceed as in Case 2, but set
$$\Sigma_1 = e_1 \cup a_1 \cup b_1 \cup \dots \cup a_g \cup b_g$$ where $e_1$ is either the sole exceptional fiber or a regular fiber.  
\end{enumerate}
\end{defn}

The following theorem is proved in \cite{MS1998}.

\begin{thm} \label{SFSHS}
Irreducible Heegaard splittings of Seifert fibered spaces are either horizontal or vertical.
\end{thm}

\section{Flipping horizontal Heegaard splittings}

\begin{thm} \label{hori}
Horizontal Heegaard splittings of Seifert fibered spaces are flippable.  
\end{thm}

\proof
Let $M = H_1 \cup_S H_2$ be a horizontal Heegaard splitting of a Seifert fibered space.  This means that there is a fiber $f$ of $M$ lying in $S$ such that $S \cap (M - \eta(f))$ is horizontal.  Moreover, as described in Remark \ref{horiHS}, $M - \eta(f)$ is a mapping torus over a surface $\hat S$, {\it i.e.,} $$M - \eta(f) =  (\hat{S} \times [-1,1])/\sim$$ and $S \cap (M - \eta(f))$ consists of two copies of $\hat S$.

The splitting surface $S$ is composed of the two copies of $\hat S$ along with a longitudinal annulus in $\eta(e)$ that is a bicollar of $e$.  The two copies of $\hat S$ are fibers of the mapping torus  $M - \eta(f) =  (\hat{S} \times [-1,1])/\sim$ and can be rotated around the mapping torus.  Furthermore, this isotopy of $M - \eta(f)$ extends into $\eta(e)$ to rotate the longitudinal annulus $S \cap \eta(f)$.  This isotopy interchanges $H_1$ and $H_2$.  The Heegaard splitting is thus flippable.  
\qed

\section{Quotients of the fundamental group of a Seifert fibered space} \label{calculation}

A Seifert fibered space $M$ can be completely described in terms of a set of invariants computed from the base orbifold, along with an invariant called the {\em Euler number}, and the $\mu$s and $\nu$s of the exceptional fibers.  For instance, when both the Seifert fibered space and the base orbifold are orientable, the invariants are as follows:

$$M = \langle g, b, (\alpha_1, \beta_1), \dots, (\alpha_l, \beta_l) \rangle$$

The number $g$ is the genus of the base orbifold $Q$ of $M$.  The number $b$ is the Euler number.  The pair $(\alpha_i, \beta_i)$ is related to $(\mu_i, \nu_i)$ but provides the slope of the meridian disk of a solid torus neighborhood of the $i$th exceptional fiber in terms of ``external" coordinates.  Specifically, denote the surface obtained by deleting a regular neighborhood of a point from $\hat Q$ by $\hat Q^-$.  The surface $\hat Q^-$ thus obtained has boundary components $(c_0, \dots, c_l)$.  For each $i$, choose a point $p_i$ in $c_i$.  Then the pair $(c_i \times \{p_i\}, \{p_i\} \times {\mathbb S}^1)$ of curves in $c_i \times {\mathbb S}^1 \subset \partial (\hat Q^- \times {\mathbb S}^1)$ provides coordinates in which we express the slope of the meridian disk of the $i$th exceptional fiber.  Take $\alpha_i$ to be the number of times a meridian of $V_i$ meets $\{p_i\} \times {\mathbb S}^1$ (in fact, $\alpha_i = \mu_i$) and $\beta_i$ to be the number of times it meets $c_i \times \{p_i\}$.  To reconstruct $M$ from $\hat Q^- \times {\mathbb S}^1$, perform Dehn fillings of slope $(\alpha_i, \beta_i)$ along $c_i \times {\mathbb S}^1$ for $i = 1, \dots, l$, then perform a Dehn filling of slope $(1, b)$ along $c_0 \times {\mathbb S}^1$.

The fundamental group of $M$ can be computed from this set of invariants.  

$$\pi_1(M) = \langle  x_1, \dots, x_l, a_1, b_1, \dots, a_g, b_g, h \; | $$
$$\; h^{-b}\Pi_1^g [a_i, b_i] \Pi_1^l x_i, [a_1, h], [b_1, h], \dots, [a_g, h], [b_g, h], $$
$$[x_1, h], \dots, [x_l, h], x_1^{\alpha_1}h^{\beta_1}, \dots, x_l^{\alpha_l}h^{\beta_l} \rangle $$

We will be interested in certain quotients of the fundamental group of $M$, for instance, the quotient by $<h>$: 

$$\pi_1(M)/<h> = \langle x_1, \dots, x_l, a_1, b_1, \dots, a_g, b_g \; | \; \Pi_1^g [a_i, b_i] \Pi_1^k x_i, x_1^{\alpha_1}, \dots, x_l^{\alpha_l} \rangle $$

In the case where there is exactly one exceptional fiber, we will also be interested in the, much simpler, quotient by the ``horizontal" subgroup $<a_1, b_1, \dots, a_g, b_g>$:

$$\pi_1(M)/<a_1, b_1, \dots, a_g, b_g> = \langle x_1, h \; | \; h^{-b}x_1, x_1^{\alpha_1}h^{\beta_1} \rangle = \langle h \; | \;  h^{b\alpha_1 + \beta_1} \rangle$$

\section{Circle bundles over surfaces}



\begin{thm} \label{bundles}
Heegaard splittings of circle bundles over a surface are flippable.
\end{thm}

\proof
Some circle bundles (those with Euler number $\pm 1$) admit horizontal Heegaard splittings.  See \cite{MS1998}.  Horizontal Heegaard splittings are flippable by Theorem \ref{SFSHS}.  It therefore suffices to consider vertical Heegaard splittings of circle bundles over surfaces.  

Let $M$ be a circle bundle over $F$.  Let $\hat F$ be a once punctured copy of $F$.  We obtain $M$ by Dehn filling $\hat F \times {\mathbb S}^1$.  Denote the core of the Dehn filling by $f_1$.  By construction, $f_1$ lies in $H_1$.  Note that $M - \eta(f_1) = (H_1 - \eta(f_1)) \cup_S H_2$ coincides with the Heegaard splitting of any other Dehn filling of $\hat F \times {\mathbb S}^1$, in particular, the Dehn filling yielding a product.  Moreover, since $f_1$ is a regular fiber, $S$ can be isotoped across $f_1$.  The rest of the proof is then identical to the proof of Theorem \ref{prodflip}.  
\qed

\section{Nielsen equivalence of generating systems}

The concept of Nielsen equivalence for generating systems of groups provides a strategy for distinguishing between Heegaard splittings.  For details, see \cite[Chapter 3]{KMS}.  

\begin{defn} 
Let $G$ be a finitely generated group with basis $[g_1, \dots, g_l]$.  The following transformations are called {\em Nielsen} moves:
\begin{enumerate}
\item Switch $x_i$ and $x_j$;
\item Replace $x_i$ with $x_i^{-1}$;
\item Replace $x_i$ with $x_i x_j$ for $i \neq j$.
\end{enumerate}
A basis obtained from $[g_1, \dots, g_n]$ by a finite sequence of these moves is said to be {\em Nielsen equivalent} to $[g_1, \dots, g_n]$.
\end{defn} 

\begin{defn} (Alternate definition)
Let $G$ be a finitely generated group with basis $[g_1, \dots, g_n]$.  The basis $[g_1', \dots, g_n']$ of $G$ is {\em Nielsen equivalent} to $[g_1, \dots, g_n]$ if it is the image of $[g_1, \dots, g_n]$ under an automorphism of $G$.
\end{defn}

To see that the two definitions of Nielsen equivalence coincide, see \cite[Chapter 3]{KMS}.  The following theorems are due to Lustig and Moriah, see \cite{LM} and \cite{LM2020}:

\begin{thm} \label{Fg}
Let $G$ be a Fuchsian group with presentation
$$G = \langle x_1, \dots, x_l, a_1, b_1, \dots, a_g, b_g \; | \; x_i^{\alpha_i}, x_1 \dots x_m \Pi_{i=1}^g [a_i, b_i] \rangle$$ with with either $l \geq 3$ or $g \geq 1$ and pairwise relatively prime exponents $\alpha_i \geq 3$. Two generating systems of $G,$ 
$${\frak A} = \{ x_1^{\nu_1}, \dots, x_{j-1}^{\nu_{j-1}}, x_{j+1}^{\nu_{j+1}}, \dots, x_l^{\nu_l}, a_1, b_1, \dots, a_g, b_g \}$$ 
$${\frak B} = \{ x_1^{\nu_1'}, \dots, x_{k-1}^{\nu_{k-1}'}, x_{k+1}^{\nu_{k+1}'}, \dots, x_l^{\nu_l'}, a_1, b_1, \dots, a_g, b_g \}$$
are Nielsen equivalent if and only if
\begin{enumerate}
\item $j = k$ and $\nu_i = \nu_i',$ or
\item $j \neq k$ and $\nu_j = 1,$ $\nu_k = 1,$ and $\nu_i = \nu_i'$ for $i = j, k.$
\end{enumerate}
\end{thm}

The calculation of $\pi_1(M)/<h>$ in Section \ref{calculation}, under the requisite hypotheses on the invariants of the exceptional fibers, yields a group satisfying the hypotheses of this theorem.  

Let $M = H_1 \cup_S H_2$ be a vertical Heegaard splitting of the Seifert fibered space $M$.  Suppose that $M$ has exceptional fibers $\{e_1, \dots, e_m\},$ $m \geq 2,$ and $e_1 \cup \dots \cup e_k$ are contained in $H_1$.  Then the fundamental group of $H_1$ induces a generating system for $\pi_1(M)$ and hence a generating system for $\pi_1(M) / \langle h \rangle$ given by $${\frak A} = \{ x_1^{\alpha_1}, \dots, x_k^{\alpha_{k}}, x_{k+2}, \dots, x_l, a_1, b_1, \dots, a_g, b_g \}$$  On the other hand, the fundamental group of $H_2$ induces a generating system for $\pi_1(M)$ and hence a generating system for $\pi_1(M) / \langle h \rangle$ given by $${\frak B} = \{ x_l^{\alpha_l}, \dots, x_{k+2}^{\alpha_{k+2}}, x_k, \dots, x_1, a_1, b_1, \dots, a_g, b_g \}$$

\begin{thm} \label{noflip}
Let $M$ be a Seifert fibered space such that  $$\pi_1(M)/<h> = \langle x_1, \dots, x_l, a_1, b_1, \dots, a_g, b_g \; | \; x_i^{\alpha_i}, x_1 \dots x_l \Pi_{i=1}^g [a_i, b_i] \rangle$$ with with either $l \geq 3$ or $g \geq 1$ and $l \geq 2$ and pairwise relatively prime exponents $\alpha_i \geq 3$.  Then irreducible vertical Heegaard splittings of $M$ are not flippable.  
\end{thm}

\proof
Let $M = H_1 \cup_S H_2$ be an irreducible vertical Heegaard splitting of the Seifert fibered space $M$.  Suppose that $M$ has exceptional fibers $\{e_1, \dots, e_n\},$ $n \geq 2,$ and $e_1 \cup \dots \cup e_k$ are contained in $H_1$.  Then the fundamental group of $H_1$ induces a generating system for $\pi_1(M)$ and hence a generating system for $\pi_1(M) / \langle h \rangle$ given by $${\frak A} = \{ x_1^{\alpha_1}, \dots, x_k^{\alpha_{k}}, x_{k+2}, \dots, x_n, a_1, b_1, \dots, a_g, b_g \}$$  On the other hand, the fundamental group of $H_2$ induces a generating system for $\pi_1(M)$ and hence a generating system for $\pi_1(M) / \langle h \rangle$ given by $${\frak B} = \{ x_n^{\alpha_n}, \dots, x_k^{\alpha_{k+2}}, x_k, \dots, x_1, a_1, b_1, \dots, a_g, b_g \}$$
If $M = H_1 \cup_S H_2$ is flippable, then the isotopy interchanging $H_1$ and $H_2$ induces an automorphism on the fundamental group.  It follows that ${\frak A}$ and ${\frak B}$ are Nielsen equivalent, but this contradicts Theorem \ref{Fg}.  Therefore the Heegaard splitting is not flippable.  
\qed

The reasoning in Theorem \ref{noflip} applies in other cases as well, even when some of the $\alpha_i = 2$.  Consider, for instance, another theorem of Lustig and Moriah:

\begin{thm} \label{Fg'}
Let $G$ be a Fuchsian group with presentation
$$G = \langle x_1, \dots, x_m, a_1, b_1, \dots, a_g, b_g \; | \; x_i^{\alpha_i}, x_1 \dots x_m \Pi_{i=1}^g [a_i, b_i] \rangle$$ and let $m$ be the number of exponents greater than $3$.  If the number  $0 < n = l - m$ of exponents equal to $2$ is even, we assume that $m \geq 5$.  If $n$ is odd, we assume that $m \geq 7$.   Two generating systems of $G,$ 
$${\frak A} = \{ x_1^{\nu_1}, \dots, x_{j-1}^{\nu_{j-1}}, x_{j+1}^{\nu_{j+1}}, \dots, x_m^{\nu_m}, a_1, b_1, \dots, a_g, b_g \}$$ 
$${\frak B} = \{ x_1^{\nu_1'}, \dots, x_{k-1}^{\nu_{k-1}'}, x_{k+1}^{\nu_{k+1}'}, \dots, x_m^{\nu_m'}, a_1, b_1, \dots, a_g, b_g \}$$
are Nielsen equivalent if and only if
\begin{enumerate}
\item $j = k$ and $\nu_i = \nu_i',$ or
\item $j \neq k$ and $\nu_j = 1,$ $\nu_k = 1,$ and $\nu_i = \nu_i'$ for $i = j, k.$
\end{enumerate}
\end{thm}

In Theorem \ref{Fg'}, exponents are allowed to coincide.  It is therefore only possible to conclude that an irreducible Heegaard splitting is not flippable if at least two distinct exceptional fibers with relatively prime exponents lie in opposite handlebodies.  Note, however, that the excluded cases are very specific.  The technique fails only for Seifert fibered spaces with an even number of exceptional fibers that pair up into duos with the same $\alpha_i$s and only for the irreducible Heegaard splittings that distribute one exceptional fiber from each such duo into one handlebody and the other into the other handlebody.  In particular, irreducible Heegaard splittings of a Seifert fibered space with an odd number of exceptional fibers with a requisite number of $\alpha_i \neq 2$ can never be flippable.  

\section{Flipping stabilized Heegaard splittings}

\begin{thm}
Stabilized Heegaard splittings of Seifert fibered spaces are flippable.
\end{thm}

\proof
Let $M = H_1 \cup_S H_2$ be a stabilized Heegaard splitting of a Seifert fibered space.  If $M = H_1 \cup_S H_2$ is a stabilization of a horizontal Heegaard splitting, then it is flippable by Theorems \ref{hori} and \ref{flipstable}.  

Following the reasoning in the proof of Theorem \ref{bundles}, it suffices to show that exceptional fibers can be moved from one side of $S$ to the other.  Note that stabilizations can be positioned where needed.  For instance, they can be positioned to occur in small $3$-balls away from the exceptional fibers.  Denote the exceptional fibers by $e_1, \dots, e_m$ and suppose that $e_1, \dots, e_l$ lie in $H_1$ and $e_{l+1}, \dots, e_m$ lie in $H_2$. 

\begin{figure}[ht]
\includegraphics[scale=.4]{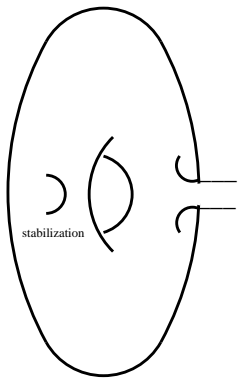}
\centering
\caption{\sl A stabilization in a small $3$-ball in a solid torus summand}
\label{flip}
\end{figure}

\begin{figure}[ht]
\includegraphics[scale=.4]{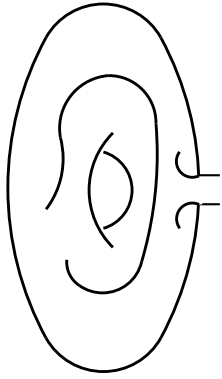}
\centering
\caption{\sl An isotopy of the core of a stabilization in a small $3$-ball in a solid torus summand}
\label{flip1}
\end{figure}

\begin{figure}[ht]
\includegraphics[scale=.4]{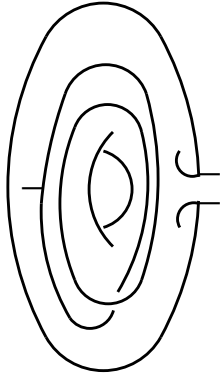}
\centering
\caption{\sl The parallel copy of the core of the solid torus in the other handlebody}
\label{flip2}
\end{figure}

We consider the exceptional fiber $e_1$ and suppose that $e_1$ lies in $H_1$.  By construction, $e_1$ is the core of a solid torus summand of $H_1$.  So $H_1 = H_1'  \cup N(a) \cup N(e_1)$, where $a$ is an arc connecting $e_1$ to the rest of $H_1$.  Next we position one stabilization to lie in a small $3$-ball inside $\partial N(e)$.  See Figure \ref{flip}.  In other words, this stabilization is realized by taking a small arc $b$ properly embedded in $N(e_1)$ and attaching a collar of $b$ to the complementary compression body.  Next we slide one of the endpoints of $b$ along $\partial N(e_1)$, see Figure \ref{flip1}, and then over the other endpoint of $b$ to obtain an eyeglass consisting of a copy of $e_1$ attached to $\partial N(e_1)$ via a short arc, see Figure \ref{flip2}.  By isotoping the simple closed curve of the eyeglass to coincide with $e_1$, we have moved $e_1$ into $H_2$.

We now (briefly) consider the Seifert fibered space with boundary $\hat M = M - \eta(e_1 \cup \dots \cup e_m)$.  We label the boundary components $B_1, \dots, B_m$ according to the deleted exceptional fiber that gave rise to them.  The splitting surface $S$ induces a Heegaard splitting of the manifold with boundary $\hat M,$ which we denote by $\hat M = \hat H_1 \cup_S \hat H_2$.  (Here $\hat H_1, \hat H_2$ are compression bodies.)  Furthermore, $B_2, \dots, B_l \subset \hat H_1$ and $B_1, B_{l+1}, \dots, B_m \subset \hat H_2$,  

The Heegaard splitting for $\hat M$ has genus larger than the genus of any unstabilized Heegaard splitting.  By \cite[Theorem 4.2]{Schultens1995}, $S$ is the splitting surface of a stabilization of a vertical Heegaard splitting of $\hat M$ and therefore also of $M$.  The partitioning of the boundary components $B_1, \dots, B_m$ guarantees that $S$ partitions $e_1, \dots, e_m$ such that $e_2, \dots, e_l  \subset H_1$ and $e_1, e_{l+1}, \dots, e_m \subset H_2$ as desired.  The exceptional fiber $e_1$ has been moved from $H_1$ into $H_2$.  
\qed


\begin{thebibliography}{BaBE}
 
\bibitem{BO1982} Bonahon, F. and Otal, J.-P., Scindement de Heegaard des espaces lenticulaires, {\em C. R. Acad. Sci. Paris S\'{e}r. I Math.} {\bf 294} (1982), no. 17, 585--587.

\bibitem{JacoBook}  Jaco, William, {\em Lectures on three-manifold topology,}
CBMS Regional Conference Series in Mathematics 43 (1980),
American Mathematical Society, Providence, R.I., xii+251 pp.
 
\bibitem{johnson2010sg} Johnson, Jesse, Bounding the stable genera of Heegaard splittings from below, {\em J. Topol.} {\bf 3} (2010), no. 3, 668--690.  

\bibitem{LM} Lustig, Martin and Moriah, Yoav, Nielsen equivalence in Fuchsian groups and Seifert fibered spaces, {\em Topology} Vol. 30, No. 2 (1991) 191--204.  

\bibitem{LM2020} Lustig, Martin and Moriah, Yoav, Nielsen equivalence in Fuchsian groups, preprint, arXiv:1910.02759v4.

\bibitem{KMS} Magnus, Wilhelm and Karrass, Abraham and Solitar, Donald, {\em Combinatorial Group Theory}, Reprint of the 1976 second edition.  Dover Publications, Inc., Minneola, NY, 2004. xii+444 pp. 

\bibitem{MS1998} Moriah, Yoav and Schultens, Jennifer, Irreducible Heegaard splittings of Seifert fibered spaces are either vertical
or horizontal,  {\em Topology} {\bf 37} (1998), no. 5, 1089--1112.

\bibitem{Schultens1993} Schultens, Jennifer, The classification of Heegaard splittings for $(compact\; orientable\; surface) \times {\mathbb S}^1$, {\em Proc. LMS} (3) {\bf 67} (1993), no, 2, 425--448.  

 \bibitem{SeifertBook} Seifert, Herbert,
 Topologie dreidimensionaler gefaserter R\"{a}ume, 
{\em Acta Math} {\bf 60} (1932) 147--238.
 
 \bibitem{HeilTrans} Seifert, Herbert and Threlfall, William,
 	{\em Seifert and Threlfall: a textbook of topology.} Translated from the German edition of 1934 by Michael A. Goldman. With a preface by Joan S. Birman. With ``Topology of 3-dimensional fibered spaces'' by Seifert. Translated from the German by Wolfgang Heil
	 {\em Pure and Applied Mathematics} {\bf 89}. Academic Press, Inc. [Harcourt Brace Jovanovich, Publishers], New York-London, 1980. xvi+437 pp.
 

\end{thebibliography}
\end{document}